\documentclass{amsart}
\usepackage[utf8]{inputenc}
\usepackage{amssymb,latexsym}
\usepackage{amsmath}
\usepackage{graphicx}
\usepackage{textcomp}

\usepackage{amsthm,amssymb,enumerate,graphicx, tikz}
\usepackage{amscd}
\usepackage{setspace}
\usepackage{comment}
\usepackage{hyperref}
\usepackage{cleveref}

\newtheorem{theorem}{Theorem}[section]

\newtheorem*{theorem*}{Theorem}

\newtheorem{lemma}[theorem]{Lemma}

\newtheorem{conjecture}[theorem]{Conjecture}
\theoremstyle{definition}

\theoremstyle{remark}

\newcommand{\F}{\mathcal{F}}

\newcommand{\conv}{\textrm{conv}}

\newcommand{\R}{\mathbb{R}}

\title{Line transversals in families of connected sets in the plane}
\author{Daniel McGinnis}\thanks{D. McGinnis: Department of Mathematics, Iowa State University, USA.  \url{dam1@iastate.edu}. D. McGinnis was supported by NSF grant DMS-1839918 (RTG)} 
\author{Shira Zerbib}\thanks{S. Zerbib: Department of Mathematics, Iowa State University, USA.  \url{zerbib@iastate.edu}. S. Zerbib was supported by NSF grant DMS-1953929.}

\begin{document}

\begin{abstract}
We prove that if a family of compact connected sets in the plane has the property that every three members of it are intersected by a line, then there are three lines intersecting all the sets in the family. This answers a question of Eckhoff from 1993 \cite{Eckhoff3}, who proved that, under the same condition, there are four lines intersecting all the sets. In fact, we prove a colorful version of this result, under weakened conditions on the sets. 
    A triple of sets $A,B,C$ in the plane is said to be a {\em tight} if $\conv(A\cup B)\cap \conv(A\cup C)\cap \conv(B\cap C)\neq \emptyset.$ This notion  was first introduced by Holmsen in \cite{Holmsen2}, where he showed that if $\F$ is a family of compact  convex sets in the plane in which every three sets form a tight triple, then there is a line intersecting at least $\frac{1}{8}|\F|$ members of $\F$. Here we prove that if $\F_1,\dots,\F_6$ are families of compact connected sets in the plane such that every three sets, chosen from three distinct families $\F_i$, form a tight triple, then there exists $1\le j\le 6$ and three lines  intersecting every member of $\F_j$. In particular, this improves $\frac{1}{8}$ to $\frac{1}{3}$ in Holmsen's result. 
\end{abstract}

\maketitle
\section{Introduction}
Let $\F$ be a family of sets in the plane. We say that $\F$ has {\em property $T(r)$} if every $r$ or fewer sets in $\F$ admit  a {\em line transversal}, that is, there exists a  line intersecting these sets. We say that $\F$ is {\em pierced} by $k$ lines if there are $k$ lines in the plane whose union intersects all the sets in $\F$. The {\em line-piercing number} of the family is the minimum $k$ so that $\F$ is pierced by $k$ lines.

The problem of bounding the line-piercing numbers of families of compact convex sets in the plane with the $T(r)$ property has been investigated since the 1960's. In 1964 Eckhoff \cite{Eckhoff1} proved that if a family of compact convex sets  satisfies the $T(4)$ property then it can be pierced by two lines. In 1974 he gave an example of a family of compact convex sets satisfying the $T(3)$ property that is not pierced by two lines \cite{Eckhoff2}. 
Various upper bounds on the line-piercing numbers were proved when further restrictions on the sets are imposed (see \cite{Eckhoff3} for more details). 

For a while it was not clear whether the $T(3)$ property implies a finite universal upper bound on the line-piercing number. This question was eventually resolved   
in 1975 by Kramer \cite{Kramer}, who showed that a family of compact convex sets in $\R^2$ with the $T(3)$ property is pierced by 5 lines. 
Finally,  in 1993 
Eckhoff \cite{Eckhoff3} proved that such families are pierced by 4 lines, and asked whether this bound can be improved to 3.  

Quantitative versions have been  studied too. In 1980, Katchalski and Liu \cite{KL} showed the existence of a constant $0<\alpha(r)<1$, so that every finite family $\F$ of  convex sets in the plane  with  the $T(r)$ property admits a line  intersecting $\alpha(r)|\F|$ of its members. In 2010  Holmsen \cite{Holmsen} showed that $\big(\frac{2}{r(r-1)}\big)^{\frac{1}{r-2}}\le \alpha(r) \le \frac{r-2}{r-1}$, and in particular,  $\frac{1}{3} \le \alpha(3) \le \frac{1}{2}.$

In \cite{Holmsen2}, Holmsen introduced the notion of tight triples. Three compact, connected sets in the plane $A,B,C$ are said to be a {\em tight triple} if $$\conv(A\cup B)\cap \conv(A\cup C)\cap \conv(B\cap C)\neq \emptyset.$$ We will call a family of sets in the plane a {\em family of tight triples} if every three sets in the family form a tight triple. Note that if $A,B,C$ has a line transversal, then it is a tight triple. Holmsen \cite{Holmsen2} proved that if $\F$ is a family of tight triples in which every set is compact and convex, then there is a line intersecting at least $\frac{1}{8}|\F|$ members of $\F$.  

The sets investigated in all the above results are assumed to be convex, but the results apply also for families of  connected sets. This follows from the fact that 
if $S$ is a  connected set in $\R^2$ and $\ell$ is a line  intersecting $\text{conv}(S)$ then $\ell$ must intersect  $S$. Similarly,  three connected sets $A,B,C$ form a tight triple if and only if $\conv(A), \conv(B), \conv(C)$ form a tight triple.

In this paper we show that the line-piercing number of a family of compact connected tight triples is at most 3. This improves the $\frac{1}{8}$ in Holmsen's result to $\frac{1}{3}$. 

In fact, we show a colorful version of this fact.

\begin{theorem}\label{main}
Let $\F_1,\dots,\F_6$ be families of compact connected sets in $\R^2$. If every three sets $A_1\in \F_{i_1}, A_2\in \F_{i_2}, A_3\in \F_{i_3}$,  $1\le i_1<i_2<i_3\le 6$, form a tight triple, then there exists  $i\in [6]$ such that the line-piercing number of $\F_i$ is at most $3$.
\end{theorem}

In particular, if $\F$ is a family with the $T(3)$ property and we take $\F_i=\F$ for $1\le i\le 6$, then Theorem \ref{main} implies that $\F$ has line-piercing number at most three. This gives an affirmative answer to Eckhoff's question.
This also gives another proof to Holmsen's result $\alpha(3) \ge 1/3$. 

Our main tool is the colorful version of the topological KKM theorem \cite{kkm} due to Gale \cite{Gale}. 
Let $\Delta^{n-1} = \{(x_1,\dots,x_n) \in \R^n \mid x_i\ge 0, \sum_{i=1}^n x_i = 1\}$ denote the $(n-1)$-dimensional simplex in $\R^n$, whose vertices are the canonical basis vectors $e_1,\dots,e_n$. Let $S_n$ be the group of permutations on $[n]$. 
\begin{theorem}[The colorful KKM theorem \cite{Gale}]\label{thm:kkm}
Let $A_1^i,\dots,A_{n}^i$, $i\in [n]$, be  open sets of $\Delta^{n-1}$, such that for every $i\in [n]$ and for every face $\sigma$ of $\Delta^{n-1}$ we have $\sigma \subset \bigcup_{e_j \in \sigma} A_j^i$.  Then there exists a permutation $\pi\in S_n$ such that $\cap_{i=1}^{n} A^{\pi(i)}_{i}\neq \emptyset$. 
\end{theorem}

\section{Proof of Theorem \ref{main}}
Throughout the proof, addition in integers is taken modulo 6. For $a,b \in \R^2$, let $[a,b]=\text{conv}\{a,b\}$ be the line segment connecting  $a,b$.

As is explained in \cite{Eckhoff3}, the compactness of the sets in each $\F_j$ allows us to assume that $\F_j$ is finite. 
Thus we may scale the plane so that every set in $\F_j$ is contained in the unit disk $D$ for each $j$.  Denote by $U$ the unit circle. Let $f(t)$ be a parameterization of $U$ defined by $f(t)=(\textrm{cos}(2\pi t), \textrm{sin}(2\pi t))$. 

A point $x=(x_1,\dots,x_6)\in \Delta^5$ corresponds to $6$ points on $U$ given by $f_i(x)=f(\sum_{j=1}^i x_{j})$ for $1\leq i\leq 6$. Let $l_1(x)=l_4(x)=[f_1(x),f_4(x)]$, $l_2(x)=l_5(x)=[f_2(x),f_5(x)]$, and $l_3(x)=l_6(x)=[f_3(x),f_6(x)]$.

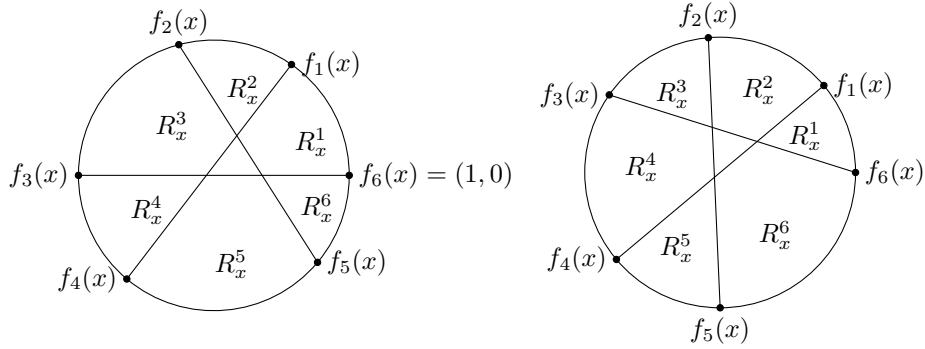
\begin{figure}[h!]
    \centering
    \begin{tikzpicture}[scale=.6,rotate=-40]
    \draw (0,0) circle (3cm);
\filldraw [black] (3,0) circle (2pt) node[right] {$f_5(x)$};
\filldraw [black] (95:3cm) circle (2pt) node[right] {$f_1(x)$};
\filldraw [black] (145:3cm) circle (2pt) node[above] {$f_2(x)$};
\filldraw [black] (270:3cm) circle (2pt) node[left] {$f_4(x)$};
\filldraw [black] (40:3cm) circle (2pt) node[right] {$f_6(x)=(1,0)$};
\filldraw [black] (220:3cm) circle (2pt) node[left] {$f_3(x)$};

\draw (3,0) -- (145:3cm);
\draw (95:3cm) -- (270:3cm);
\draw (40:3cm) -- (220:3cm);

\filldraw [] (325:1.5cm) node[below] {$R^5_x$};
\filldraw [] (35:2.3cm) node[below] {$R^6_x$};
\filldraw [] (115:2.5cm) node[below] {$R^2_x$};
\filldraw [] (160:1.85cm) node[below] {$R^3_x$};
\filldraw [] (70:2.5cm) node[below] {$R^1_x$};
\filldraw [] (230:1.5cm) node[below] {$R^4_x$};
\filldraw [] (310:3cm) node[below] {$\phantom{f} $};

\end{tikzpicture}
\begin{tikzpicture}[scale=.6]
    \draw (0,0) circle (3cm);
\filldraw [black] (3,0) circle (2pt) node[right] {$f_6(x)$};
\filldraw [black] (95:3cm) circle (2pt) node[above] {$f_2(x)$};
\filldraw [black] (145:3cm) circle (2pt) node[left] {$f_3(x)$};
\filldraw [black] (270:3cm) circle (2pt) node[below] {$f_5(x)$};
\filldraw [black] (40:3cm) circle (2pt) node[right] {$f_1(x)$};
\filldraw [black] (220:3cm) circle (2pt) node[left] {$f_4(x)$};

\draw (3,0) -- (145:3cm);
\draw (95:3cm) -- (270:3cm);
\draw (40:3cm) -- (220:3cm);

\filldraw [] (325:1.5cm) node[below] {$R^6_x$};
\filldraw [] (35:2.3cm) node[below] {$R^1_x$};
\filldraw [] (115:2.5cm) node[below] {$R^3_x$};
\filldraw [] (160:1.85cm) node[below] {$R^4_x$};
\filldraw [] (70:2.5cm) node[below] {$R^2_x$};
\filldraw [] (230:1.5cm) node[below] {$R^5_x$};

\end{tikzpicture}
 \caption{A point $x\in \Delta^5$ corresponds to six regions $R_x^i$. The regions $R^1_p,R^3_p,R^5_p$ are pairwise disjoint (on the right) or the regions $R^2_p,R^4_p,R^6_p$ are pairwise disjoint (on the left), depending on the orientation of the triangle bounded by the lines $l_1(x),l_2(x),l_3(x)$.}
    \label{figure1}
\end{figure}

For $i=1,\dots, 6$ let $R^i_x$ be the interior of the region bounded by  $l_{i-1}(x)$, $l_{i}(x)$ and the arc on $U$ connecting $f_{i-1}(x)$ and $f_{i}(x)$ (see Figure \ref{figure1}).  Notice that $R^i_x=\emptyset$ when $x_{i}=0$. Also, it is possible that some of the regions $R^i_x$ intersect.

Set $1\le j \le 6$ and let $A_i^j$  be the set of points $x\in \Delta^5$ so that $R^i_x$ contains a set $F\in \F_j$. Since the sets $F\in \F_j$ are closed, $A_i^j$ is open.  
If there is some $x\in \Delta^5$ for which $x\notin \bigcup_{i=1}^6 A_i^j$, then since the sets in $\F_j$ are connected, every set in $\F_j$ must intersect $\bigcup_{i=1}^3 l_i(x)$, and we are done. 
So we assume for contradiction that $\Delta^5=\bigcup_{i=1}^6 A_i^j$ for all $j$. Observe that if 
 $x\in \textrm{conv}\{e_i:i\in I\}$ for some $I\subset [6]$ then $R^k_x=\emptyset$ for $k\notin I$, and  therefore, $x \in \bigcup_{i\in I} A_i^j$ for all $j$. This shows that the conditions of Theorem \ref{thm:kkm} hold.

Thus, by Theorem \ref{thm:kkm}, there exists some permutation $\pi \in S_6$ and a point $p = (p_1, \dots, p_6) \in \bigcap_{i=1}^6A_{i}^{\pi(i)}$. Therefore, each of the open regions $R^i_p$ contains a set $S_i \in \F_{\pi(i)}$, $i = 1, \dots, 6$, and in particular $R^i_p \neq \emptyset$ and thus $p_i\neq 0$ for all $i$. We claim that at least one of the triples $\{S_1, S_3, S_5\}$ or $\{S_2,S_4,S_6\}$ is not a tight triple. To see this,  
note that the regions $R^1_p,R^3_p,R^5_p$ are pairwise disjoint or the regions $R^2_p,R^4_p,R^6_p$ are pairwise disjoint (depending on the orientation of the triangle bounded by the lines $l_1,l_2,l_3$, see Figure \ref{figure1}). Without loss of generality, we assume $R^1_p,R^3_p,R^5_p$ are pairwise disjoint, and in this case, the three sets $S_1, S_3, S_5$  is not a tight triple. This is a contradiction.

 

 


\section{Concluding remarks}

The proof of Theorem \ref{main} implies a slightly stronger result: when each $\F_i$ is finite, one can fix a point lying on one of the three piercing lines of $\F_{i}$, as long as this point is outside $\text{conv}(\cup_i\F_i)$.

A similar proof can be used to prove a colorful version of Eckhoff's result that $T(4)$ families are pierced by two lines. \begin{theorem}
Let $\F_1,\dots,\F_4$ be families of compact, connected sets in the plane such that any collection of four sets, one from each $\F_i$,  has a line transversal. Then for some $i\in [4]$, $\F_{i}$ has line piercing number at most 2.
\end{theorem}
This can be proved by associating a point in $\Delta^3$ with two lines and applying a similar argument as in the proof of Theorem \ref{main}

When each of the families $\F_i$ is finite, one may drop the condition that the sets are compact. This is because we may replace each set $S\in \F_j$ with a compact, convex set $S'\subset \text{conv}(S)$ such that the resulting family is still a family of tight triples. 

\section{Acknowledgement}
We are grateful to Andreas Holmsen for many  helpful discussions and for telling us about Eckhoff's paper \cite{Eckhoff3}. We are also grateful to Ron Aharoni for commenting on an early version of this paper.

\end{document}